\documentclass[12pt,a4paper]{article}
\usepackage[cp1251]{inputenc}
\usepackage{amsfonts, amssymb, amsthm, amsmath}
\usepackage{mathtext}
\usepackage{longtable}
\usepackage[russian]{babel}
\usepackage{color}
\usepackage{hyperref}
\hypersetup{
colorlinks=true,
linkcolor=blue,
citecolor=blue,
urlcolor=black
}

\usepackage{cite}
\usepackage[a4paper, margin=2cm]{geometry}

\begin{document}

\renewcommand{\r}{\mathbb R}
\newcommand{\nb}{N\!\!\!B}
\newcommand{\N}{\mathbb N}
\newcommand{\eqd}{\stackrel{d}{=}}
\newcommand{\pto}{\stackrel{P}{\longrightarrow}}

\title{On mixture representations for the generalized Linnik distribution
and their applications in limit theorems}

\author{V. Yu. Korolev\textsuperscript{1}, A. K.
Gorshenin\textsuperscript{2}, A. I. Zeifman\textsuperscript{3}}

\date{}

\maketitle

\footnotetext[1]{Hangzhou Dianzi University, China;
Faculty of Computational Mathematics and Cybernetics, Lomonosov
Moscow State University; Institute of Informatics Problems, Federal
Research Center ``Computer Science and Control'' of the Russian
Academy of Sciences; vkorolev@cs.msu.ru}

\footnotetext[2]{Institute of Informatics Problems, Federal Research
Center ``Computer Science and Control'' of the Russian Academy of
Sciences; Faculty of Computational Mathematics and Cybernetics,
Lomonosov Moscow State University; agorshenin@frccsc.ru}

\footnotetext[3]{Vologda State University; Institute of Informatics Problems, Federal
Research Center ``Computer Science and Control'' of the Russian
Academy of Sciences; a$\_$zeifman@mail.ru}

{

\small

{\bf Abstract:} We present new mixture representations for the
generalized Linnik distribution in terms of normal, Laplace,
exponential and stable laws and establish the relationship between
the mixing distributions in these representations. Based on these
representations, we prove some limit theorems for a wide class of
rather simple statistics constructed from samples with random sized
including, e. g., random sums of independent random variables with
finite variances and maximum random sums, in which the generalized
Linnik distribution plays the role of the limit law. Thus we
demonstrate that the scheme of geometric (or, in general, negative
binomial) summation is far not the only asymptotic setting (even for
sums of independent random variables) in which the generalized
Linnik law appears as the limit distribution.

\smallskip


}

\section{Introduction}

In this paper we continue the research we started in
\cite{KorolevZeifman2016, KorolevZeifmanKMJ}. We study the
interrelationship between the (generalized) Linnik and (generalized)
Mittag-Leffler distributions. In \cite{KorolevZeifmanKMJ} we showed
that along with the traditional and well-known representation of the
Linnik distribution as the scale mixture of a strictly stable law
with exponential mixing distribution, there exists another
representation of the Linnik law as the normal scale mixture with
the Mittag-Leffler mixing distribution. The former representation
makes it possible to treat the Linnik law as the limit distribution
for geometric random sums of independent identically distributed
random variables (r.v.'s) in which summands have infinite variances.
The latter normal scale mixture representation opens the way to
treating the Linnik distribution as the limit distribution in the
central limit theorem for random sums of independent random
variables in which summands have {\it finite} variances. Moreover,
being scale mixtures of normal laws, the Linnik distributions can
serve as the one-dimensional distributions of a special subordinated
Wiener process often used as models of the evolution of stock prices
and financial indexes. Strange as it may seem, the results
concerning the possibility of representation of the Linnik
distribution as a scale mixture of normals were never explicitly
presented in the literature in full detail before the paper
\cite{KorolevZeifmanKMJ} saw the light, although the property of the
Linnik distribution to be a normal scale mixture is something almost
obvious. Perhaps, the paper \cite{KotzOstrovskii1996} was the
closest to this conclusion and exposed the representability of the
Linnik law as a scale mixture of Laplace distributions with the
mixing distribution written out explicitly. Methodically, the
present paper is very close to the work of L. Devroye
\cite{Devroye1996} where many examples of mixture representations of
popular probability distributions were discussed from the simulation
point of view.

Here we consider the generalized Linnik distribution and prove new
mixture representations for it in terms of normal, Laplace,
exponential and stable laws and establish the relationship between
the mixing distributions in these representations. In particular, we
prove that the generalized Linnik distribution is a normal scale
mixture with the generalized Mittag-Leffler mixing distribution.
Based on these representations, we prove some limit theorems for a
wide class of rather simple statistics constructed from samples with
random sized including, e. g., random sums of independent random
variables with finite variances and maximum random sums, in which
the generalized Linnik distribution plays the role of the limit law.
Thus we demonstrate that the scheme of geometric (or, in general,
negative binomial) summation is far not the only asymptotic setting
(even for sums of independent r.v.'s) in which the generalized
Linnik law appears as the limit distribution.

The presented material substantially relies on the results of
\cite{KorolevZeifmanKMJ} and \cite{LimTeo2009}. The paper is
organized as follows. Section 2 presents the definitions and basic
properties of the Linnik, generalized Linnk, and Mittag-Leffler
distributions. Section 3 contains basic definitions and auxiliary
results. In Section 4 we prove the representability of the
generalized Linnik distribution as the normal scale mixture of
normal laws with the generalized Mittag-Leffler mixing distribution.
We show that if the `generalizing' parameter $\nu$ does not exceed
1, then the generalized Linnik distribution is a scale mixture of
`ordinary' Linnik distributions with the same characteristic
parameter $\alpha$. We prove the representation of the generalized
Linnik distribution as a scale mixture of the Laplace laws with the
mixing distribution explicitly determined as that of the randomly
scaled ratio of two independent r.v.'s with the same strictly stable
distribution concentrated on the nonnegative halfline. Here we also
discuss some properties of the generalized Mittag-Leffler
distribution. In Sections 5 and 6 we prove and discuss some criteria
(that is, necessary and sufficient conditions) for the convergence
of the distributions of rather simple statistics constructed from
samples with random sizes including, e. g., random sums of
independent r.v.'s with finite variances to the generalized Linnik
law.

\section{The Linnik and Mittag-Leffler distributions}

\subsection{The Linnik distributions}

In this paper our attention will be focused on the distributions
whose characteristic functions (ch.f.) have the form
$$
\mathfrak{f}^L_{\alpha}(t)=\frac{1}{1+|t|^{\alpha}},\ \ \
t\in\mathbb{R},\eqno(1)
$$
where $\alpha\in(0,2]$. These symmetric distributions were
introduced 1n 1953 by Yu. V. Linnik \cite{Linnik1953}. The
distributions with the ch.f. (1) are traditionally called the {\it
Linnik distributions}. Although sometimes the term {\it
$\alpha$-Laplace distributions} \cite{Pillai1985} is used, we will
use the first term which has already become conventional. If
$\alpha=2$, then the Linnik distribution turns into the Laplace
distribution corresponding to the density
$$
f^{\Lambda}(x)=\textstyle{\frac12}e^{-|x|},\ \ \
x\in\mathbb{R}.\eqno(2)
$$
A r.v. with density (2) and the corresponding distribution function
(d.f.) will be denoted $\Lambda$ and $F^{\Lambda}(x)$, respectively.
A r.v. with the Linnik distribution with parameter $\alpha$ will be
denoted $L_{\alpha}$. Its d.f. will be denoted $F_{\alpha}^{L}$.
From (1) and (2) we obviously have $F_2^{L}(x)\equiv
F^{\Lambda}(x)$, $x\in\mathbb{R}$.

The properties of the Linnik distributions were studied in many
papers. We should mention \cite{Laha1961, Devroye1990,
KotzOstrovskiiHayfavi1995a, KotzOstrovskiiHayfavi1995b} and other
papers, see the survey in \cite{KorolevZeifmanKMJ}.

Perhaps, most often Linnik distributions are recalled as examples of
geometric stable distributions. This means that if $X_1,X_2,\ldots$
are independent r.v.'s whose distributions belong to the domain of
attraction of an $\alpha$-strictly stable symmetric law and
$\nb_{1,\,p}$ is the r.v. independent of $X_1,X_2,\ldots$ and having
the geometric distribution
$$
{\sf P}(\nb_{1,\,p}=n)=p(1-p)^{n-1},\ \ \ n=1,2,\ldots,\ \ \
p\in(0,1),
$$
then for each $p\in(0,1)$ there exists a constant $a_p>0$ such that
$a_p\big(X_1+\ldots+X_{\nb_{1,\,p}}\big)\Longrightarrow L_{\alpha}$
as $p\to 0$, see, e. g., \cite{Bunge1996} or
\cite{KlebanovRachev1996} (the symbol $\Longrightarrow$ hereinafter
denotes convergence in distribution).

In \cite{Pakes1998}, Pakes showed that the probability distributions
known as {\it generalized Linnik distributions} which have
characteristic functions
$$
\mathfrak{f}^L_{\alpha,\nu,\theta}(t)=\frac{1}{(1+e^{-i\theta\,{\mathrm{sgn}\,t}}|t|^{\alpha})^{\nu}},\
\ \ t\in\mathbb{R},\
|\theta|\le\min\Big\{\frac{\pi\alpha}{2},\,\pi-\frac{\pi\alpha}{2}\Big\},\
\nu>0,\eqno(3)
$$
play an important role in some characterization problems of
mathematical statistics. The class of probability distributions with
density functions $p_{\alpha,\nu,\theta}(x)$ corresponding to ch.f.
(3) have found some interesting properties and applications, see
\cite{Anderson1993, BaringhausGrubel1997, Devroye1990,
Jayakumar1995, KotzOstrovskii1996, Kotz2001, Kozubowski1998,
Lim1998} and related papers. In particular, they are good candidates
to model financial data which exhibits high kurtosis and heavy tails
\cite{MittnikRachev1993}.

Here we concentrate our attention at the symmetric case $\theta=0$.
Any r.v. with the ch.f.
$$
\mathfrak{f}^L_{\alpha,\nu,0}(t)=\frac{1}{(1+|t|^{\alpha})^{\nu}},\
\ \ t\in\mathbb{R},
$$
will be denoted $L_{\alpha,\nu}$.

\subsection{The Mittag-Leffler distributions}

The Mittag-Leffler probability distribution is the distribution of a
nonnegative r.v. $M_{\delta}$ whose Laplace transform is
$$
\psi_{\delta}(s)\equiv {\sf E}e^{-sM_{\delta}}=\frac{1}{1+\lambda
s^{\delta}},\ \ \ s\ge0,\eqno(4)
$$
where $\lambda>0$, $0<\delta\le1$. For simplicity, in what follows
we will consider the standard scale case and assume that
$\lambda=1$.

The origin of the term {\it Mittag-Leffler distribution} is due to
that the probability density corresponding to Laplace transform (4)
has the form
$$
f_{\delta}^{M}(x)=\frac{1}{x^{1-\delta}}\sum\nolimits_{n=0}^{\infty}\frac{(-1)^nx^{\delta
n}}{\Gamma(\delta n+1)}=-\frac{d}{dx}E_{\delta}(-x^{\delta}),\ \ \
x\ge0,
$$
where $E_{\delta}(z)$ is the Mittag-Leffler function with index
$\delta$ that is defined as the power series
$$
E_{\delta}(z)=\sum\nolimits_{n=0}^{\infty}\frac{z^n}{\Gamma(\delta
n+1)},\ \ \ \delta>0,\ z\in\mathbb{Z}.
$$
Here $\Gamma(s)$ is Euler's gamma-function,
$$
\Gamma(s)=\int_{0}^{\infty}z^{s-1}e^{-z}dz,\ \ \ s>0.
$$
The Mittag-Leffler distribution function (d.f.) will be denoted
$F_{\delta}^{M}(x)$, $F_{\delta}^{M}(x)={\sf P}(M_{\delta}<x)$,
$x\in\mathbb{R}$.

With $\delta=1$, the Mittag-Leffler distribution turns into the
standard exponential distribution, that is, $F_1^{M}(x)=
[1-e^{-x}]\mathbf{1}(x\ge 0)$, $x\in\mathbb{R}$ (here and in what
follows the symbol $\mathbf{1}(C)$ denotes the indicator function of
a set $C$). But with $\delta<1$ the Mittag-Leffler distribution
density has the heavy power-type tail: from the well-known
asymptotic properties of the Mittag-Leffler function it can be
deduced that if $0<\delta<1$, then
$$
f_\delta^{M}(x)\sim \frac{\sin(\delta\pi)\Gamma(\delta+1)}{\pi
x^{\delta+1}}
$$
as $x\to\infty$, see, e. g., \cite{Kilbas2014}.

It is well-known that the Mittag-Leffler distribution is
geometrically stable. The history of the Mittag-Leffler distribution
is discussed in \cite{KorolevZeifmanKMJ}.

The Mittag-Leffler distributions are of serious theoretical interest
in the problems related to thinned (or rarefied) homogeneous flows
of events such as renewal processes or anomalous diffusion or
relaxation phenomena, see \cite{WeronKotulski1996,
GorenfloMainardi2006} and the references therein.

Let $\nu>0$, $\delta\in(0,1]$. By analogy to (3), the distribution
of a nonnegative r.v. $M_{\delta,\,\nu}$ defined by the
Laplace--Stieltjes transform
$$
\psi_{\delta,\,\nu}(s)\equiv {\sf
E}e^{-sM_{\delta,\,\nu}}=\frac{1}{(1+s^{\delta})^{\nu}},\ \ \ s\ge0,
$$
will be called the {\it generalized Mittag-Leffler distribution},
see \cite{MathaiHaubold2011, Joseetal} and the references therein.
Sometimes this distribution is called the {\it Pillai distribution}
\cite{Devroye1996}, although in the original paper \cite{Pillai1985}
R. Pillai called it {\it semi-Laplace}. In the present paper we will
keep to the first term {\it generalized Mittag-Leffler
distribution}.

\section{Basic notation and auxiliary results}

Most results presented below actually concern special mixture
representations for probability distributions. However, without any
loss of generality, for the sake of visuality and compactness of
formulations and proofs we will represent the results in terms of
the corresponding r.v.'s assuming that all the r.v.'s mentioned in
what follows are defined on the same probability space
$(\Omega,\,\mathfrak{A},\,{\sf P})$.

The r.v. with the standard normal d.f. $\Phi(x)$ will be denoted
$X$,
$$
{\sf
P}(X<x)=\Phi(x)=\frac{1}{\sqrt{2\pi}}\int_{-\infty}^{x}e^{-z^2/2}dz,\
\ \ \ x\in\mathbb{R}.
$$
Let $\Psi(x)$, $x\in\mathbb{R}$, be the d.f. of the maximum of the
standard Wiener process on the unit interval,
$\Psi(x)=2\Phi\big(\max\{0,x\}\big)-1$, $x\in\mathbb{R}$. It is easy
to see that $\Psi(x)={\sf P}(|X|<x)$. Therefore, sometimes $\Psi(x)$
is said to determine the {\it half-normal} or {\it folded normal}
distribution.

A r.v. having the gamma distribution with shape parameter $r>0$ and
scale parameter $\lambda>0$ will be denoted $G_{r,\lambda}$,
$$
{\sf P}(G_{r,\lambda}<x)=\int_{0}^{x}g(z;r,\lambda)dz,\ \
\text{with}\ \
g(x;r,\lambda)=\frac{\lambda^r}{\Gamma(r)}x^{r-1}e^{-\lambda x},\
x\ge0,
$$
where $\Gamma(r)$ is Euler's gamma-function,
$\Gamma(r)=\int_{0}^{\infty}x^{r-1}e^{-x}dx$, $r>0$.

In this notation, obviously, $G_{1,1}$ is a r.v. with the standard
exponential distribution: ${\sf P}(G_{1,1}<x)=\big[1-e^{-x}\big]{\bf
1}(x\ge0)$ (here and in what follows ${\bf 1}(A)$ is the indicator
function of a set $A$).

The gamma distribution is a particular representative of the class
of generalized gamma distributions (GG distributions), which were
first described in \cite{Stacy1962} as a special family of lifetime
distributions containing both gamma and Weibull distributions.

A {\it generalized gamma $($GG$)$ distribution} is the absolutely
continuous distribution defined by the density
$$
\overline{g}(x;r,\alpha,\lambda)=\frac{|\alpha|\lambda^r}{\Gamma(r)}x^{\alpha
r-1}e^{-\lambda x^{\alpha}},\ \ \ \ x\ge0,
$$
with $\alpha\in\mathbb{R}$, $\lambda>0$, $r>0$.

A r.v. with the density $\overline{g}(x;r,\alpha,\lambda)$ will be
denoted $\overline{G}_{r,\alpha,\lambda}$.

It is easy to see that
$$
\overline{G}_{r,\alpha,\mu}\eqd G_{r,\mu}^{1/\alpha}\eqd
\mu^{-1/\alpha}G_{r,1}^{1/\alpha}\eqd\mu^{-1/\alpha}\overline{G}_{r,\alpha,1}.\eqno(5)
$$

The d.f. and the density of the strictly stable distribution with
the characteristic exponent $\alpha$ and shape parameter $\theta$
defined by the characteristic function
$$
\mathfrak{p}_{\alpha,\theta}(t)=\exp\big\{-|t|^{\alpha}\exp\{-{\textstyle\frac12}i\pi\theta\alpha\,\mathrm{sign}t\}\big\},\
\ \ \ t\in\r,\eqno(6) 
$$
with $0<\alpha\le2$, $|\theta|\le\min\{1,\frac{2}{\alpha}-1\}$, will
be denoted by $P_{\alpha,\theta}(x)$ and $p_{\alpha,\theta}(x)$,
respectively (see, e. g., \cite{Zolotarev1983}). Any r.v. with the
d.f. $P_{\alpha,\theta}(x)$ will be denoted $S_{\alpha,\theta}$. For
definiteness, $S_{1,\,1}\eqd 1$ (throughout the paper the symbol
$\eqd$ will denote the coincidence of distributions).

From (6) it follows that the characteristic function of a symmetric
($\theta=0$) strictly stable distribution has the form
$$
\mathfrak{p}_{\alpha,0}(t)=e^{-|t|^{\alpha}},\ \ \ t\in\r. \eqno(7) 
$$
From (7) it is easy to see that $S_{2,0}\eqd\sqrt{2}X$.

\smallskip

{\sc Lemma 1}. (i) {\it Let $\alpha\in(0,2]$, $\alpha'\in(0,1]$.
Then}
$$
S_{\alpha\alpha',0}\eqd S_{\alpha,0}S_{\alpha',1}^{1/\alpha}
$$
{\it where the r.v.'s on the right-hand side are independent.}

\noindent (ii) {\it Let $\alpha\in(0,1]$, $\alpha'\in(0,1]$. Then}
$$
S_{\alpha\alpha',1}\eqd S_{\alpha,1}S_{\alpha',1}^{1/\alpha}
$$
{\it where the r.v.'s on the right-hand side are independent.}

\smallskip

{\sc Proof}. See, e. g., Theorem 3.3.1 in \cite{Zolotarev1983}.

\smallskip

{\sc Corollary 1.} {\it A symmetric strictly stable distribution
with the characteristic exponent $\alpha$ is a scale mixture of
normal laws in which the mixing distribution is the one-sided
strictly stable law $(\theta=1)$ with the characteristic exponent
$\alpha/2:$
$$
S_{\alpha,0}\eqd X\sqrt{2S_{\alpha/2,1}} \eqno(8)
$$
with the r.v.'s on the right-hand side being independent.}

\smallskip

In terms of d.f.s the statement of Corollary 1 can be written as
$$
P_{\alpha,0}(x)=\int_{0}^{\infty}\Phi\Big(\frac{x}{\sqrt{2z}}\Big)dP_{\alpha/2,1}(z),\
\ \ x\in\r.
$$

\smallskip

Let $\gamma>0$. The distribution of the r.v. $W_{\gamma}$:
$$
{\sf
P}\big(W_{\gamma}<x\big)=\big[1-e^{-x^{\gamma}}\big]\mathbf{1}(x\ge
0),
$$
is called the {\it Weibull distribution} with shape parameter
$\gamma$. It is obvious that $W_1$ is the r.v. with the standard
exponential distribution: ${\sf P}(W_1<x)=\big[1-e^{-x}\big]{\bf
1}(x\ge0)$. The Weibull distribution is a particular case of GG
distributions corresponding to the density
$\overline{g}(x;1,\gamma,1)$. Thus, $W_1\eqd G_{1,1}$. It is easy to
see that $W_1^{1/\gamma}\eqd W_{\gamma}$.

It is easy to see that if $\gamma>0$ and $\gamma'>0$, then ${\sf
P}(W_{\gamma'}^{1/\gamma}\ge x)={\sf P}(W_{\gamma'}\ge
x^{\gamma})=e^{-x^{\gamma\gamma'}}={\sf P}(W_{\gamma\gamma'}\ge x)$,
$x\ge 0$, that is, for any $\gamma>0$ and $\gamma'>0$
$$
W_{\gamma\gamma'}\eqd W_{\gamma'}^{1/\gamma}.\eqno(9) 
$$
In particular, $W_{\gamma}\eqd W_1^{1/\gamma}$.

Any Weibull distribution with $0<\gamma\le1$ is a mixed exponential
distribution. Namely, if $0<\gamma\le1$, then
$$
W_{\gamma}\eqd\frac{W_1}{S_{\gamma,1}},\eqno(10)
$$
where the r.v.'s on the right-hand side are independent (see Lemma 4
in \cite{Korolev2016} or Lemma 3 in \cite{KorolevZeifmanKMJ}).

\smallskip

In the paper \cite{Gleser1989} it was shown that any gamma
distribution with shape parameter no greater than one is mixed
exponential. For convenience, we formulate this result as the
following lemma.

\smallskip

{\sc Lemma 2} \cite{Gleser1989}. {\it The density $g(x;r,\mu)$ of a
gamma distribution with $0<r<1$ can be represented as
$$
g(x;r,\mu)=\int_{0}^{\infty}ze^{-zx}p(z;r,\mu)dz,
$$
where
$$
p(z;r,\mu)=\frac{\mu^r}{\Gamma(1-r)\Gamma(r)}\cdot\frac{\mathbf{1}(z\ge\mu)}{(z-\mu)^rz}.\eqno(11)
$$
Moreover, a gamma distribution with shape parameter $r>1$ cannot be
represented as a mixed exponential distribution.}

\smallskip

{\sc Lemma 3} \cite{Korolev2017}. {\it For $r\in(0,1)$ let
$G_{r,\,1}$ and $G_{1-r,\,1}$ be independent gamma-distributed
r.v.'s. Let $\mu>0$. Then the density $p(z;r,\mu)$ defined by $(11)$
corresponds to the r.v.
$$
Z_{r,\mu}=\frac{\mu(G_{r,\,1}+G_{1-r,\,1})}{G_{r,\,1}}\eqd\mu
Z_{r,1}\eqd\mu\big(1+\textstyle{\frac{1-r}{r}}V_{1-r,r}\big),
$$
where $V_{1-r,r}$ is the r.v. with the Snedecor--Fisher distribution
defined by the probability density
$$
q(x;1-r,r)=\frac{(1-r)^{1-r}r^r}{\Gamma(1-r)\Gamma(r)}
\cdot\frac{1}{x^{r}[r+(1-r)x]},\ \ \ x\ge0.
$$
}

\smallskip

{\sc Remark 1}. It is easily seen that the sum $G_{r,1}+G_{1-r,1}$
has the standard exponential distribution: $G_{r,1}+G_{1-r,1}\eqd
W_1$. However, the numerator and denominator of the expression
defining the r.v. $Z_{r,\mu}$ are not independent.

\smallskip

The statements of Lemmas 2 and 3 mean that if $r\in(0,1)$, then
$$
G_{r,\,\mu}\eqd W_1\cdot Z_{r,\,\mu}^{-1}\eqno(12)
$$
where the r.v.'s on the righ-hand side are independent.

In \cite{Korolev2017} (also see \cite{KorolevZeifman2019LMJ}) the
following statement similar to Lemma 2 was proved for generalized
gamma distributions.

\smallskip

{\sc Lemma 4} \cite{Korolev2017, KorolevZeifman2019LMJ}. {\it Let
$\alpha\in(0,1]$, $r\in(0,1)$, $\mu>0$. Then the GG distribution
with parameters $r$, $\alpha$, $\mu$ is a mixed exponential
distribution$:$
$$
\overline{G}_{r,\alpha,\mu}\eqd
W_1\cdot\big(S_{\alpha,1}Z_{r,\mu}^{1/\alpha}\big)^{-1}
$$
with the r.v.'s on the right-hand side being independent. Moreover,
a GG distribution with $\alpha r>1$ cannot be represented as mixed
exponential.}

\smallskip

The following statement has already become folklore. Without any
claims for priority, its proof was given in \cite{KorolevZeifmanKMJ}
as an exercise.

\smallskip

{\sc Lemma 5.} {\it For any $\delta\in(0,1]$, the Mittag-Leffler
distribution with parameter $\delta$ is a scale mixture of a
one-sided stable distribution with the Weibull mixing distribution
with parameter $\delta$, that is,}
$$
M_{\delta}\eqd S_{\delta,1}W_{\delta},
$$
{\it where the r.v.'s on the right-hand side are independent}.

\smallskip

From Lemma 5 and (10) we obtain
$$
M_{\delta}\eqd S_{\delta,1}\cdot W_{\delta}\eqd
W_1\cdot\frac{S_{\delta,1}}{S'_{\delta,1}},
$$
where all r.v.'s on the right-hand side are independent. Denote
$$
R_{\delta}=\frac{S_{\delta,1}}{S'_{\delta,1}}.
$$
In \cite{KorolevZeifmanKMJ} it was shown that the probability
density $f^{R}_{\delta}(x)$ of the ratio $R_{\delta}$ of two
independent r.v.'s with one and the same one-sided strictly stable
distribution with parameter $\delta$ has the form
$$
f^{R}_{\delta}(x)=
\frac{\sin(\pi\delta)x^{\delta-1}}{\pi[1+x^{2\delta}+2x^{\delta}\cos(\pi\delta)]},\
\ \ x>0.\eqno(13)
$$
So, the following statement is valid.

\smallskip

{\sc Lemma 6} \cite{Kozubowski1998, KorolevZeifmanKMJ}. {\it Let
$0<\delta\le1$. The Mittag-Leffler distribution is mixed
exponential$:$
$$
M_{\delta}\eqd W_1\cdot R_{\delta},
$$
where the r.v.'s on the right-hand side are independent and the
probability density of the mixing distribution is given by $(13)$}.

\smallskip

The Mittag-Leffler distribution with $\delta<1$ is a scale mixture
of the Mittag-Leffler distributions with a greater parameter, as is
seen from the following statement.

\smallskip

{\sc Lemma 7}. {\it Let $0<\delta\le1$, $0<\delta'\le1$. Then
$$
M_{\delta\delta'}\eqd M_{\delta}\cdot R_{\delta'}^{1/\delta},
$$
where $R_{\delta'}$ is the ratio of two independent r.v.'s with one
and the same one-sided strictly stable distribution with parameter
$\delta'$ and the r.v.'s on the right-hand side are independent}.

{\sc Proof}. From Lemmas 5 and 1(ii), relations (9) and (10) we
obtain the following chain of relations:
$$
M_{\delta\delta'}\eqd S_{\delta\delta',1}\cdot
W_1^{1/\delta\delta'}\eqd S_{\delta,1}\cdot
S_{\delta',1}^{1/\delta}\cdot W_1^{1/\delta\delta'}\eqd
S_{\delta,1}\cdot S_{\delta',1}^{1/\delta}\cdot
W_{\delta'}^{1/\delta}\eqd
$$
$$
\eqd S_{\delta,1}\cdot
S_{\delta',1}^{1/\delta}\cdot\Big(\frac{W_1}{S'_{\delta',1}}\Big)^{1/\delta}\eqd
S_{\delta,1}\cdot W_1^{1/\delta}\cdot R_{\delta'}^{1/\delta}\eqd
M_{\delta}\cdot R_{\delta'}^{1/\delta}.
$$
The lemma is proved.

\smallskip

It is easily seen that the probability density of the mixing
distribution of the r.v. $R_{\delta'}^{1/\delta}$ has the form
$\delta x^{\delta-1}f^R_{\delta'}(x^{\delta})$ with $f^R$ given by
(13).

Lemma 7 is a concretization of a result of \cite{Kozubowski1998}
where it was demonstrated that the Mittag-Leffler distribution with
$\delta<1$ is a scale mixture of the Mittag-Leffler distributions
with a greater parameter and the explicit form of the mixing density
was presented, but the mixing distribution was not recognized as the
distribution of the powered ratio of two independent random
variables with one and the same one-sided strictly stable
distribution.

Now turn to the Linnik distribution. In \cite{Devroye1990} the
following statement was proved. Here its formulation is extended
with the account of (9).

\smallskip

{\sc Lemma 8} \cite{Devroye1990}. {\it For any $\alpha\in(0,2]$, the
Linnik distribution with parameter $\alpha$ is a scale mixture of a
symmetric stable distribution, that is,}
$$
L_{\alpha}\eqd S_{\alpha,0}\cdot W_1^{1/\alpha},
$$
{\it where the r.v.'s on the right-hand side are independent}.

\smallskip

In the same way as Lemma 7 was proved, it can be shown that any
Linnik distribution with $0<\alpha<2$ is a scale mixture of the
Linnik distributions with a greater parameter.

Let $0<\alpha\le2$, $0<\alpha'\le1$. Then from Lemmas 8 and 1(i),
relations (9) and (10) we obtain the following chain of relations:
$$
L_{\alpha\alpha'}\eqd S_{\alpha\alpha',0}\cdot
W_1^{1/\alpha\alpha'}\eqd S_{\alpha,0}\cdot
S_{\alpha',1}^{1/\alpha}\cdot W_1^{1/\alpha\alpha'}\eqd
S_{\alpha,0}\cdot S_{\alpha',1}^{1/\alpha}\cdot
W_{\alpha'}^{1/\alpha}\eqd
$$
$$
\eqd S_{\alpha,0}\cdot
S_{\alpha',1}^{1/\alpha}\cdot\Big(\frac{W_1}{S'_{\alpha',1}}\Big)^{1/\alpha}\eqd
S_{\alpha,0}\cdot W_1^{1/\alpha}\cdot R_{\alpha'}^{1/\alpha}\eqd
L_{\alpha}\cdot R_{\alpha'}^{1/\alpha}.
$$
Therefore, the following statement holds.

\smallskip

{\sc Lemma 9}. {\it Let $0<\alpha\le2$, $0<\alpha'\le1$. Then
$$
L_{\alpha\alpha'}\eqd L_{\alpha}\cdot R_{\alpha'}^{1/\alpha},
$$
where $R_{\alpha'}$ is the ratio of two independent r.v.'s with one
and the same one-sided strictly stable distribution with parameter
$\alpha'$ and the r.v.'s on the right-hand side are independent}.

\smallskip

Lemma 9 is a concretization of a result of \cite{KotzOstrovskii1996}
where it was demonstrated that the Linnik distribution with
$\alpha<2$ is a scale mixture of the Linnik distributions with a
greater parameter and the explicit form of the mixing density was
presented, but the mixing distribution was not recognized as the
distribution of the powered ratio of two independent random
variables with one and the same one-sided strictly stable
distribution.

\smallskip

With $\alpha=2$ from Lemma 9 we obtain

\smallskip

{\sc Corollary 3} \cite{KorolevZeifmanKMJ}. {\it Let $0<\alpha<2$.
Then the Linnik distribution with parameter $\alpha$ is a scale
mixture of the Laplace distributions corresponding to density
$(2)$}:
$$
L_{\alpha}\eqd \Lambda \cdot \sqrt{R_{\alpha/2}},
$$
{\it where the r.v.'s on the right-hand side are independent.}

\smallskip

Now we recall some representations of the Linnik distribution as a
normal or Laplace scale mixtures from \cite{KorolevZeifmanKMJ}.

\smallskip

{\sc Lemma 10} \cite{KorolevZeifmanKMJ}. {\it Let $\alpha\in(0,2]$,
$\alpha'\in(0,1]$. Then}
$$
L_{\alpha\alpha'}\eqd S_{\alpha,0}\cdot M_{\alpha'}^{1/\alpha}
$$
{\it where the r.v.'s on the right-hand side are independent}.
\smallskip

As far as we know, although the property of the Linnik distribution
to be a normal scale mixture is something almost obvious by virtue
of Lemmas 8 and 1, only in the paper \cite{KorolevZeifmanKMJ} the
mixing distribution was recognized as the Mittag-Leffler law, as is
seen from the following statement.

\smallskip

{\sc Corollary 4} \cite{KorolevZeifmanKMJ}. {\it For each
$\alpha\in(0,2]$, the Linnik distribution with parameter $\alpha$ is
the scale mixture of zero-mean normal laws with mixing
Mittag-Leffler distribution with twice less parameter $\alpha/2$}:
$$
L_{\alpha}\eqd X\sqrt{2M_{\alpha/2}},\eqno(14)
$$
{\it where the r.v.'s on the right-hand side are independent}.

\smallskip

Now consider the folded normal mixture representation for the
Mittag-Leffler distribution.

\smallskip

{\sc Lemma 11} \cite{KorolevZeifmanKMJ}. {\it For $\delta\in(0,1]$
the Mittag-Leffler distribution with parameter $\delta$ is a scale
mixture of half-normal laws}:
$$
M_{\delta}\eqd \sqrt{2}|X|\cdot R_{\delta}\cdot W_2.
$$

\smallskip

Notice that the statement of Lemma 11 can also be interpreted as
that the Mittag-Leffler distribution is a scale mixture of the
Rayleigh laws.

\section{New mixture representations for the generalized Linnik
distribution}

In this section we present some results containing new mixture
representations for the generalized Linnik distribution. These
results generalize and improve some results of \cite{Pakes1998} and
\cite{LimTeo2009}. We begin from the following well-known result due
to Devroye and Pakes.

\smallskip

{\sc Lemma 12} \cite{Devroye1990, Pakes1998}. {\it Let
$\alpha\in(0,2]$, $\nu>0$. Then}
$$
L_{\alpha,\nu}\eqd S_{\alpha,0}\cdot G_{\nu,1}^{1/\alpha}\eqd
S_{\alpha,0}\cdot\overline{G}_{\nu,\alpha,1}.
$$

\smallskip

{\sc Remark 2}. Let $D_{\nu}$ be a r.v. with the {\it one-sided
exponential power distribution} defined by the density
$$
f^D_{\nu}(x)=\frac{1}{\nu\Gamma(1/\nu)}e^{-x^{1/\nu}},\ \ \ x\ge0.
$$
It is easy to make sure that $D_{\nu}^{1/\nu}\eqd G_{\nu,1}$.
Therefore, the statement of Lemma 12 can be re-formulated as
$$
L_{\alpha,\nu}\eqd S_{\alpha,0}\cdot D_{\nu}^{1/\alpha\nu},
$$
as it was done in \cite{Devroye1996}.

\smallskip

From Corollary 1 (see relation (8)) it follows that for $\nu>0$ and
$\alpha\in(0,2]$
$$
L_{\alpha,\nu}\eqd X\cdot\sqrt{2S_{\alpha/2,1}}\cdot
G_{\nu,1}^{1/\alpha}\eqd
X\cdot\sqrt{2S_{\alpha/2,1}\cdot\overline{G}_{\nu,\alpha/2,1}}.\eqno(15)
$$
that is, the generalized Linnik distribution is a normal scale
mixture.

Notice that the Laplace--Stieltjes transform of the mixing
distribution in (15) has the form
$$
\psi(s;\nu,\alpha/2)={\sf
E}e^{-sS_{\alpha/2,1}\overline{G}_{\nu,\,\alpha/2,\,1}}={\sf E}{\sf
E}\big[e^{-sS_{\alpha/2,1}\overline{G}_{\nu,\,\alpha/2,\,1}}\big|\overline{G}_{\nu,\,\alpha/2,\,1}\big]=
$$
$$
=\int_{0}^{\infty}{\sf
E}e^{-xsS_{\alpha/2,1}}\overline{g}(x;\nu,\alpha/2,1)dx=
\frac{\alpha}{2\Gamma(\nu)}\int_{0}^{\infty}e^{-x^{\alpha/2}(1+s^{\alpha/2})}x^{\alpha
\nu/2-1}dx=
$$
$$
=\frac{\alpha}{2\Gamma(\nu)(1+s^{\alpha/2})^\nu}
\int_{0}^{\infty}e^{-x^{\alpha/2}}x^{\alpha
\nu/2-1}dx=\frac{1}{(1+s^{\alpha/2})^\nu},\ \ \ s\ge0,
$$
corresponding to the generalized Mittag-Leffler distribution with
parameters $\alpha/2$ and $\nu$, that is,
$$
M_{\alpha/2,\,\nu}\eqd
S_{\alpha/2,1}\overline{G}_{\nu,\,\alpha/2,\,1}\eqno(16)
$$
(see \cite{MathaiHaubold2011, Joseetal}). So, by analogy with
Corollary 4 we obtain the following result.

\smallskip

{\sc Theorem 1.} {\it If $\alpha\in(0,2]$ and $\nu>0$, then
$$
L_{\alpha,\nu}\eqd X\cdot\sqrt{2M_{\alpha/2,\,\nu}},
$$
where the r.v.'s on the right-hand side are independent. In other
words, the generalized Linnik distribution is a normal scale mixture
with the generalized Mittag-Leffler mixing distribution.}

\smallskip

It is easy to see that for any $\alpha>0$ and $\alpha'>0$
$$
\overline{G}_{\nu,\alpha\alpha',1}\eqd
G_{\nu,1}^{1/\alpha\alpha'}\eqd
(G_{\nu,1}^{1/\alpha'})^{1/\alpha}\eqd\overline{G}_{\nu,\alpha',1}^{1/\alpha}.
$$
Therefore, for $\alpha\in(0,2]$, $\alpha'\in(0,1)$ and $\nu>0$ using
Lemma 1 and (16) we obtain the following chain of relations:
$$
L_{\alpha\alpha',\,\nu}\eqd S_{\alpha\alpha',\,0}\cdot
G_{\nu,1}^{1/\alpha\alpha'}\eqd S_{\alpha,\,0}\cdot
S_{\alpha',\,1}^{1/\alpha}\cdot G_{\nu,1}^{1/\alpha\alpha'}\eqd
S_{\alpha,\,0}\cdot
(S_{\alpha',\,1}\overline{G}_{\nu,\alpha',1})^{1/\alpha}\eqd
S_{\alpha,0}\cdot M_{\alpha',\nu}^{1/\alpha}.
$$
Hence, the following statement, more general than Theorem 1, holds
representing the generalized Linnik distribution as a scale mixture
of a symmetric stable law with any greater characteristic parameter,
the mixing distribution being the generalized Mittag-Leffler law.

\smallskip

{\sc Theorem 2.} {\it Let $\alpha\in(0,2]$, $\alpha'\in(0,1)$ and
$\nu>0$. Then}
$$
L_{\alpha\alpha',\,\nu}\eqd S_{\alpha,0}\cdot
M_{\alpha',\nu}^{1/\alpha}.
$$

\smallskip

Now let $\nu\in(0,1]$. From representation (12) and Lemma 8 we
obtain the chain
$$
L_{\alpha,\nu}\eqd S_{\alpha,0}\cdot G_{\nu,1}^{1/\alpha}\eqd
S_{\alpha,0}\cdot W_1^{1/\alpha}\cdot Z_{\nu,1}^{-1/\alpha}\eqd
S_{\alpha,0}\cdot W_{\alpha}\cdot Z_{\nu,1}^{-1/\alpha}\eqd
L_{\alpha}\cdot Z_{\nu,1}^{-1/\alpha}
$$
yielding the following statement relating {\it generalized} and
`ordinary' Linnik distributions.

\smallskip

{\sc Theorem 3}. {\it If $\nu\in(0,1]$ and $\alpha\in(0,2]$, then}
$$
L_{\alpha,\nu}\eqd L_{\alpha}\cdot Z_{\nu,1}^{-1/\alpha}.\eqno(17)
$$

\smallskip

In other words, with $\nu\in(0,1]$ and $\alpha\in(0,2]$, the
generalized Linnik distribution is a scale mixture of ordinary
Linnik distributions.

\smallskip

{\sc Remark 3}. Unlike the case of `ordinary' Linnik laws (see Lemma
9), the representation of the generalized Linnik law as a scale
mixture of generalized Linnik law with a greater parameter is not so
transparent. Namely, in \cite{Pakes1998} it was proved that if
$0<\alpha\le2$, $0<\alpha'<1$ and $\nu>0$, then
$$
L_{\alpha\alpha',\nu}\eqd L_{\alpha,\nu}\cdot
S_{\alpha',1}^{1/\alpha}\cdot V_{\alpha,\alpha',\nu},
$$
where
$$
V_{\alpha,\alpha',\nu}\eqd \exp\{-K_{\alpha,\alpha',\nu}\}
$$
and the distribution of the r.v. $K_{\alpha,\alpha',\nu}$ is defined
by its Laplace--Stiltjes transform
$$
\psi^{(K)}_{\alpha,\alpha',\nu}(s)=\exp\bigg\{-\int_{0}^{\infty}(1-e^{-su})\pi_{\alpha,\alpha',\nu}(u)du\bigg\}
$$
with
$$
\pi_{\alpha,\alpha',\nu}(u)=\frac{1}{u}\bigg(\frac{e^{-\nu\alpha\alpha'u}}{1-e^{-\alpha\alpha'u}}-
\frac{e^{-\nu\alpha u}}{1-e^{-\alpha u}}\bigg).
$$

\smallskip

From (17) and Lemma 10 we obtain the following representation of the
generalized Linnik distribution as a scale mixture of Laplace
distributions.
$$
L_{\alpha,\nu}\eqd \Lambda\cdot
Z_{\nu,1}^{-1/\alpha}\cdot\sqrt{R_{\alpha/2}}.
$$

Furthermore, from Corollary 4 (see relation (14)) it follows that,
if $\nu\in(0,1)$ and $\alpha\in(0,2]$, then
$$
L_{\alpha,\nu}\eqd X\cdot
Z_{\nu,1}^{-1/\alpha}\cdot\sqrt{2M_{\alpha/2}}.\eqno(18)
$$

Since normal scale mixtures are identifiable \cite{Teicher1961},
from representation (18) and Theorem 1 we obtain the following
result representing the {\it generalized} Mittag-Leffler
distribution as a scale mixture of `ordinary' Mittag-Leffler
distributions.

\smallskip

{\sc Theorem 4}. {\it Let $\nu\in(0,1]$ and $\delta\in(0,1]$. Then
$$
M_{\delta,\,\nu}\eqd Z_{\nu,1}^{-1/\delta}\cdot M_{\delta},
$$
where the r.v.'s on the right-hand side are independent}.

\smallskip

Let $\alpha\in(0,2]$, $\alpha'\in(0,1]$. From Theorems 1, 2 and
Corollary 1 we obtain the following chain of relations:
$$
X\cdot\sqrt{2M_{\alpha\alpha'/2,\nu}}\eqd
L_{\alpha\alpha'/2,\nu}\eqd S_{\alpha,0}\cdot
M_{\alpha',\nu}^{1/\alpha}\eqd X\cdot\sqrt{2S_{\alpha/2,1}}\cdot
M_{\alpha',\nu}^{1/\alpha}.
$$
Now replacing here $\alpha/2$ by $\delta\in(0,1]$ and $\alpha'$ by
$\delta'\in(0,1]$ we obtain the relation
$$
X\cdot\sqrt{2M_{\delta\delta',\nu}}\eqd
X\cdot\sqrt{2S_{\delta,1}}\cdot M_{\delta',\nu}^{1/2\delta}.
$$
Hence by virtue of the identifiability of scale mixtures of normal
laws we obtain the following statement.

\smallskip

{\sc Theorem 5.} {\it Let $\delta\in(0,1]$, $\delta'\in(0,1]$ and
$\nu>0$. Then
$$
M_{\delta\delta',\,\nu}\eqd S_{\delta,1}\cdot
M_{\delta',\nu}^{1/\delta},
$$
where the r.v.'s on the right-hand side are independent}.

\smallskip

That is, any generalized Mittag-Leffler distribution is a scale
mixture of a one-sided stable law with any greater characteristic
parameter, the mixing distribution being the generalized
Mittag-Leffler law.

\section{Convergence of the distributions of random sums to the
generalized Linnik distribution}

In applied probability it is a convention that a model distribution
can be regarded as well-justified or adequate, if it is an {\it
asymptotic approximation}, that is, if there exists a rather simple
limit setting (say, schemes of maximum or summation of random
variables) and the corresponding limit theorem in which the model
under consideration manifests itself as a limit distribution. The
existence of such limit setting can provide a better understanding
of real mechanisms that generate observed statistical regularities.

As we have already mentioned, the Linnik distributions are
geometrically stable. Geometrically stable distributions are only
possible limits for the distributions of geometric random sums of
independent identically distributed r.v.'s. As this is so, the
distributions of the summands belong to the domain of attraction of
the strictly stable law with some characteristic exponent
$\alpha\in(0,2]$ and hence, for $0<\alpha<2$ have infinite moments
of orders greater or equal to $\alpha$. As concerns the case
$\alpha=2$, where the variance is finite, within the framework of
the scheme of geometric summation in this case the only possible
limit law is the Laplace distribution \cite{Klebanov}.

As we will demonstrate below, the generalized Linnik distributions
can be limiting for negative binomial sums of independent
identically distributed r.v.'s. Negative binomial random sums turn
out to be important and adequate models of total precipitation
volume during wet (rainy) periods in meteorology \cite{Gulev,
KorolevGorsheninDAN}. However, in this case the summands also must
have distributions from the domain of attraction of a strictly
stable law with some characteristic exponent $\alpha\in(0,2]$ and
hence, with $\alpha\in(0,2)$, have infinite variances. If
$\alpha=2$, then the only possible limit distribution for negative
binomial random sums is the so-called variance gamma distribution
which is well known in financial mathematics
\cite{GnedenkoKorolev1996}.

However, when the (generalized) Linnik distributions are used as
models of statistical regularities observed in real practice and an
additive structure model is used of type of a (stopped) random walk
for the observed process, the researcher cannot avoid thinking over
the following question: which of the two combinations of conditions
can be encountered more often:

\begin{itemize}

\item the distribution of the number of summands (the number of jumps of a random walk)
is negative binomial (asymptotically gamma), but the distributions
of summands (jumps) have so heavy tails that, at least, their
variances are infinite, or

\item the second moments (variances) of the summands (jumps) are finite, but the
number of summands exposes an irregular behavior so that its very
large values are possible?

\end{itemize}

Since, as a rule, when real processes are modeled, there are no
serious reasons to reject the assumption that the variances of jumps
are finite, the second combination at least deserves a thorough
analysis.

As it was demonstrated in the preceding section, the (generalized)
Linnik distributions even with $\alpha<2$ can be represented as
normal scale mixtures. This means that they can be limit
distributions in analogs of the central limit theorem for random
sums of independent r.v.'s {\it with finite variances}. Such analogs
with `ordinary' Linnik limit distributions were presented in
\cite{KorolevZeifmanKMJ}. Here we will extend these results to
generalized Linnik distributions. It will de demonstrated that the
scheme of negative binomial summation is far not the only asymptotic
setting (even for sums of independent r.v.'s!) in which the
generalized Linnik law appears as a limit distribution.

We will begin with the limit theorem for negative binomial random
sums in which the generalized Linnik law is the limit distribution.
For this purpose we will use the following auxiliary result.
Consider a sequence of r.v.'s $Y_1,Y_2,...$ Let $N_1,N_2,...$ be
natural-valued r.v.'s such that for every $n\in\mathbb{N}$ the r.v.
$N_n$ is independent of the sequence $Y_1,Y_2,...$ In the following
statement the convergence is meant as $n\to\infty$. The symbol
$\Longrightarrow$ will denote convergence in distribution. Recall
that a random sequence $N_1,N_2,\ldots$ is said to infinitely
increase in probability ($N_n\pto\infty$), if ${\sf P}(N_n\le
m)\longrightarrow 0$ for any $m\in(0,\infty)$.

\smallskip

{\sc Lemma 13} \cite{Korolev1994, Korolev1995}. {\it Assume that
there exist an infinitely increasing $($convergent to zero$)$
sequence of positive numbers $\{b_n\}_{n\ge1}$ and a r.v. $Y$ such
that
$$
b_n^{-1}Y_n\Longrightarrow Y.\eqno(19)
$$
If there exist an infinitely increasing $($convergent to zero$)$
sequence of positive numbers $\{d_n\}_{n\ge1}$ and a r.v. $N$ such
that
$$
d_n^{-1}b_{N_n}\Longrightarrow N,\eqno(20)
$$
then
$$
d_n^{-1}Y_{N_n}\Longrightarrow Y\cdot N,\eqno(21)
$$
where the r.v.'s on the right-hand side of $(21)$ are independent.
If, in addition, $N_n\longrightarrow\infty$ in probability and the
family of scale mixtures of the d.f. of the r.v. $Y$ is
identifiable, then condition $(20)$ is not only sufficient for
$(21)$, but is necessary as well.}

\smallskip

Let $X_1,X_2,\ldots$ be independent identically distributed random
variables such that their common distribution belongs to the domain
of attraction of a strictly stable law with the characteristic
exponent $\alpha\in(0,2]$. This means that there exists a
$c\in(0,\infty)$ such that
$$
\frac{1}{cn^{1/\alpha}}\sum\nolimits_{i=1}^nX_i\Longrightarrow
S_{\alpha,\,0}\eqno(22)
$$
as $n\to\infty$. For simplicity, without any restriction of
generality, we will assume that $c=1$.

Consider a r.v. $\nb_{\nu,p}$ having the negative binomial
distribution with parameters $\nu>0$ and $p\in(0,1)$:
$$
{\sf
P}(\nb_{\nu,p}=k)=\frac{\Gamma(\nu+k-1)}{(k-1)!\Gamma(\nu)}\cdot
p^{\nu}(1-p)^{k-1},\ \ \ \ k=1,2,...,
$$
In this case ${\sf E}\nb_{\nu,p}=\nu/p$.

\smallskip

{\sc Lemma 14.} {\it If $p\to0$, then}
$$
p\nb_{\nu,p}\Longrightarrow G_{\nu,1}.
$$

\smallskip

The proof is a simple exercise on characteristic functions.

\smallskip

Let $n\in\mathbb{N}$. Assume that $p=1/n$ and denote
$N_n=\nb_{\nu,\,1/n}$. Also assume that for each $n\ge1$ the random
variable $N_n$ is independent of $X_1,X_2,\ldots$ Let
$Y_n=X_1+\ldots+X_n$. Put $b_n=d_n=n^{1/\alpha}$. Then from Lemma 14
it follows that
$$
b_{N_n}/d_n=\big(N_n/n\big)^{1/\alpha}\Longrightarrow
G_{\nu,1}^{1/\alpha}
$$
as $n\to\infty$. We will treat this relation as condition (20) in
Lemma 13. As condition (19) we will treat relation (22) (with
$c=1$). As a result, from Lemmas 12 and 14 and we obtain the
following statement establishing the convergence of the
distributions of negative binomial random sums to the generalized
Linnik distribution.

\smallskip

{\sc Theorem 6.} {\it Let $X_1,X_2,\ldots$ be independent
identically distributed r.v.'s such that their common distribution
belongs to the domain of attraction of a strictly stable law with
the characteristic exponent $\alpha\in(0,2]$. Let $N_n$ be a random
variable having the negative binomial distribution with parameters
$\nu>0$ and $p=1/n$. Then
$$
\frac{1}{n^{1/\alpha}}\sum\nolimits_{i=1}^{N_n}X_i\Longrightarrow
L_{\alpha,\,\nu}
$$
as $n\to\infty$.}

\smallskip

{\sc Remark 4.} If $\alpha=2$, then the limit r.v. has the form
$$
L_{2,\,\nu}\eqd X\cdot\sqrt{2G_{\nu,1}}.
$$
The distribution of this r.v. is a normal scale mixture with respect
to the gamma distribution. This is the well-known {\it variance
gamma distribution}, a popular heavy-tailed model in financial
mathematics.

\smallskip

As we have already noted, the representation for the generalized
Linnik distribution as a scale mixture of normals obtained above
opens the way for the construction in this section of a random-sum
central limit theorem with the generalized Linnik distribution as
the limit law. Moreover, in this ``if and only if'' version of the
random-sum central limit theorem the generalized Mittag-Leffler
distribution {\it must} be the limit law for the normalized number
of summands.

Consider independent not necessarily identically distributed random
variables $X_1,X_2,\ldots $ with ${\sf E}X_i=0$ and
$0<\sigma^2_i={\sf Var}X_i<\infty$, $i\ge1$. For $n\in\N$ denote
$$
S^*_n=X_1+\ldots +X_n,\ \ \ \  B^2_n=\sigma^2_1+\ldots +\sigma^2_n.
$$
Assume that the r.v.'s $X_1,X_2,\ldots $ satisfy the Lindeberg
condition: for any $\tau>0$
$$
\lim_{n\to\infty}\frac{1}{B^2_n}\sum\nolimits_{i=1}^{n}\int_{|x|\ge\tau
B_n}^{} x^2d{\sf P}(X_i<x)=0.\eqno(23)
$$
It is well known that under these assumptions
$$
{\sf P}\big(S^*_n<B_nx\big)\Longrightarrow \Phi(x)
$$
(this is the classical Lindeberg central limit theorem).

Let $N_1,N_2,\ldots$ be a sequence of integer-valued nonnegative
r.v.'s defined on the same probability space so that for each
$n\in\N$ the r.v. $N_n$ is independent of the sequence
$X_1,X_2,\ldots$ Denote $S^*_{N_n}=X_1+\ldots +X_{N_n}$. For
definiteness, in what follows we assume that
$\sum\nolimits_{j=1}^0=0$. Everywhere in what follows the
convergence will be meant as $n\to\infty$.

Let $\{d_n\}_{n\ge1}$ be an infinitely increasing sequence of
positive numbers.

The proof of the main result of this section is based on the
following version of the random-sum central limit theorem.

\smallskip

{\sc Lemma 15} \cite{Korolev1994}. {\it Assume that the random
variables $X_1,X_2,\ldots$ and $N_1,N_2,\ldots$ satisfy the
conditions specified above. In particular, let Lindeberg condition
$(23)$ hold. Moreover, let $N_n\pto\infty$. A d.f. $F(x)$ such that
$$
{\sf P}\Big(\frac{S^*_{N_n}}{d_n}<x\Big) \Longrightarrow F(x)
$$
exists if and only if there exists a d.f. $H(x)$ satisfying the
conditions
$$
H(0)=0,\ \ \
F(x)=\int_{0}^{\infty}\Phi\Big(\frac{x}{\sqrt{y}}\Big)dH(y),\ \
x\in\mathbb{R},
$$
and ${\sf P}(B^2_{N_n}<xd_n^2)\Longrightarrow H(x)$. }

\smallskip

{\sc Proof}. This statement is a particular case of a result proved
in \cite{Korolev1994}, also see Theorem 3.3.2 in
\cite{GnedenkoKorolev1996}.

\smallskip

The following theorem gives a {\it criterion} (that is, {\it
necessary and sufficient} conditions) of the convergence of the
distributions of random sums of independent identically distributed
r.v.'s with {\it finite} variances to the generalized Linnik
distribution.

\smallskip

{\sc Theorem 7.} {\it Let $\alpha\in(0,2]$, $\nu>0$. Assume that the
r.v.'s $X_1,X_2,\ldots$ and $N_1,N_2,\ldots$ satisfy the conditions
specified above. In particular, let Lindeberg condition $(23)$ hold.
Moreover, let $N_n\pto\infty$. Then the distributions of the
normalized random sums $S^*_{N_n}$ converge to the generalized
Linnik law with parameters $\alpha$ and $\nu$, that is,
$$
\frac{S^*_{N_n}}{d_n}\Longrightarrow L_{\alpha,\,\nu}
$$
with some $d_n>0$, $d_n\to\infty$, if and only if
$$
\frac{B^2_{N_n}}{d^2_n}\Longrightarrow 2M_{\alpha/2,\,\nu}.
$$
}

\smallskip

{\sc Proof}. This statement is a direct consequence of Theorem 1 and
Lemma 15 with $ H(x)={\sf P}(2M_{\alpha/2,\,\nu}<x)$.

\smallskip

Note that if the r.v.'s $X_1,X_2,\ldots$ are identically
distributed, then $\sigma_i=\sigma$, $i\in\N$, and the Lindeberg
condition holds automatically. In this case it is reasonable to take
$d_n=\sigma\sqrt{n}$. Hence, from Theorem 5 in this case it follows
that for the convergence
$$
\frac{S^*_{N_n}}{\sigma\sqrt{n}}\Longrightarrow L_{\alpha,\,\nu}
$$
to hold it is necessary and sufficient that
$$
\frac{N_n}{n}\Longrightarrow 2M_{\alpha/2,\,\nu}.
$$

One more remark is that with $\alpha=2$ Theorem 5 involves the case
of convergence to the Laplace distribution.

\section{Convergence of the distributions of statistics
constructed from samples with random sizes to the generalized Linnik
distribution}

In classical problems of mathematical statistics, the size of the
available sample, i. e., the number of available observations, is
traditionally assumed to be deterministic. In the asymptotic
settings it plays the role of infinitely increasing {\it known}
parameter. At the same time, in practice very often the data to be
analyzed is collected or registered during a certain period of time
and the flow of informative events each of which brings a next
observation forms a random point process. Therefore, the number of
available observations is unknown till the end of the process of
their registration and also must be treated as a (random)
observation. In this case the number of available observations as
well as the observations themselves are unknown beforehand and
should be treated as random to avoid underestimation of risks or
error probabilities.

Therefore it is quite reasonable to study the asymptotic behavior of
general statistics constructed from samples with random sizes for
the purpose of construction of suitable and reasonable asymptotic
approximations. As this is so, to obtain non-trivial asymptotic
distributions in limit theorems of probability theory and
mathematical statistics, an appropriate centering and normalization
of r.v.'s and vectors under consideration must be used. It should be
especially noted that to obtain reasonable approximation to the
distribution of the basic statistics, both centering and normalizing
values should be non-random. Otherwise the approximate distribution
becomes random itself and, for example, the problem of evaluation of
quantiles or significance levels becomes senseless.

In asymptotic settings, statistics constructed from samples with
random sizes are special cases of random sequences with random
indices. The randomness of indices usually leads to that the limit
distributions for the corresponding random sequences are
heavy-tailed even in the situations where the distributions of
non-randomly indexed random sequences are asymptotically normal see,
e. g., \cite{BeningKorolev2002, BeningKorolev2005,
GnedenkoKorolev1996}.

Consider a problem setting that is traditional for mathematical
statistics. Let r.v.'s $N_1,N_2,\ldots,X_1,X_2,\ldots,$ be defined
on one and the same probability space $(\Omega,{\cal A}, {\sf P})$.
Assume that for each $n\ge 1$ the r.v. $N_n$ takes only natural
values and is independent of the sequence $X_1,X_2,\ldots$ Let
$T_n=T_n(X_1,\ldots,X_n)$ be a statistic, that is, a measurable
function of $X_1,\ldots,X_n$. For every $n\ge1$ define the random
variable $T_{N_n}$ as
$$
T_{N_n}(\omega)=
T_{N_n(\omega)}\left(X_1(\omega),\ldots,X_{N_n(\omega)}(\omega)\right)
$$
for each $\omega\in\Omega$. As usual, the symbol $\Longrightarrow$
denotes convergence in distribution.

A statistic $T_n$ is said to be {\it asymptotically normal}, if
there exist $\sigma>0$ and $\theta\in\r$ such that
$$
{\sf
P}\left(\sigma\sqrt{n}\bigl(T_n-\theta\bigr)<x\right)\Longrightarrow\Phi(x).\eqno(24)
$$

\smallskip

{\sc Lemma 16} \cite{Korolev1995}. {\it Assume that
$N_n\longrightarrow\infty$ in probability. Let the statistic $T_n$
be asymptotically normal in the sense of $(24)$. A distribution
function $F(x)$ such that
$$
{\sf P}\left(\sigma\sqrt{n}\bigl(T_{N_n}-\theta\bigr)<x\right)
\Longrightarrow F(x),
$$
exists if and only if there exists a d.f. $H(x)$ satisfying the
conditions
$$
H(0)=0,\ \ \ F(x)=\int_{0}^{\infty}\Phi\big(x\sqrt{y}\big)dH(y),\ \
x\in\mathbb{R},\ \ \ {\sf P}(N_n<nx)\Longrightarrow H(x).
$$
}

\smallskip

The following theorem gives a criterion (that is, {\it necessary and
sufficient} conditions) of the convergence of the distributions of
statistics, which are suggested to be asymptotically normal in the
traditional sense but are constructed from samples with random
sizes, to the generalized Linnik distribution.

\smallskip

{\sc Theorem 8.} {\it Let $\alpha\in(0,2]$, $\nu>0$. Assume that the
r.v.'s $X_1,X_2,\ldots$ and $N_1,N_2,\ldots$ satisfy the conditions
specified above and, moreover, let $N_n\pto\infty$. Let the
statistic $T_n$ be asymptotically normal in the sense of $(24)$.
Then the distribution of the statistic $T_{N_n}$ constructed from
samples with random sizes $N_n$ converges to the generalized Linnik
law, that is,
$$
\sigma\sqrt{n}\bigl(T_{N_n}-\theta\bigr) \Longrightarrow
L_{\alpha,\,\nu},
$$
if and only if
$$
\frac{N_n}{n}\Longrightarrow \textstyle{
\frac12}M_{\alpha/2,\,\nu}^{-1}.
$$
}

\smallskip

{\sc Proof}. This statement is a direct consequence of Theorem 1 and
Lemma 16 with $H(x)={\sf P}(M_{\alpha/2,\,\nu}^{-1}<2x)$.

\bigskip

\subsection*{Acknowledgement} V. Yu. Korolev and A. K. Gorshenin were
supported by the Russian Foundation for Basic Research, project
17-07-00717.

\end{document}